\documentclass[12pt,letterpaper]{amsart}
\usepackage{amssymb,latexsym, amsmath, amsxtra}
\usepackage[dvips]{graphics}

\theoremstyle{plain}
\newtheorem{theorem}{Theorem}[section]

\newtheorem{lemma}[theorem]{Lemma}
\newtheorem*{rprinciple}{Repulsion Principle}
\newtheorem{proposition}[theorem]{Proposition}
\theoremstyle{definition}

\theoremstyle{remark}

\numberwithin{equation}{section}

\newcommand{\Z}{\mathbb Z}
\newcommand{\C}{\mathbb C}

\def\({\left(}
\def\){\right)}

\begin{document}
\title{Differentiating polynomials, and $\zeta(2)$}
\author{David W. Farmer \and Robert Rhoades}

\thanks{Research of the first author supported by the
American Institute of Mathematics and
the NSF Focused Research Group grant DMS 0244660}

\thispagestyle{empty} \vspace{.5cm}
\begin{abstract}
We study the derivatives of polynomials with equally spaced
zeros and find connections to the values of the
Riemann zeta-function at the integers.
\end{abstract}

\maketitle
\section{Introduction}

Polynomials are fascinating because they have many facets to their
personalities.  By definition, a polynomial $f\in \C[x]$ is an
expression of the form
\begin{equation}
f(x)=a_0+a_1 x + \cdots + a_n x^n ,
\end{equation}
where the $a_j$ are complex numbers.
By the fundamental theorem of algebra, $f$ also has a
representation as a product
\begin{equation}
f(x)=a_n (x-x_1) \cdots (x-x_n),
\end{equation}
where $x_j$ are the roots of~$f$.
The fact that there are two ways of looking at polynomials
provides possibilities that are hidden by their apparent simplicity.

In this paper we look at two different polynomials which have equally
spaced zeros,
and we study the zeros of their derivatives. We will see that
Euler's famous result
\begin{equation}\label{eqn:zeta(2)}
\zeta(2):=\sum_{n=1}^\infty \frac{1}{n^2} = \frac{\pi^2}{6}
\end{equation}
is lurking in the background, as are some generalizations of Euler's
result.  Our proof that $\zeta(2)=\pi^2/6$ appears to be new,
and it doesn't use anything more than is typically
covered in a second semester calculus course.

\section{A Polynomial With Zeros at the
Integers}\label{sec:zerosatIntegers} The first polynomial we will
consider is the degree $N+1$ polynomial with simple zeros at
$0,1,2,\dots, N$. That is,
\begin{equation}\label{eqn:product}
p_N(x) := x(x-1)\cdots (x-N) = \prod_{n=0}^N (x-n).
\end{equation}
We will write $p(x)$ instead of $p_N(x)$ to simplify the notation.
The question we will obtain a partial answer to is:
Where are the zeros of the derivative $p'(x)$?

By Rolle's theorem, we know that $p'$ has a zero in each of the intervals
$(0,1)$, $(1,2)$,\,\ldots,$(N-1,N)$.  Those zeros account for $N$ of the zeros
of $p'$, which is all of the zeros because $p'$ has degree~$N$.

We are left with the question:  where in those intervals
$(0,1)$,\ldots,$(N-1,N)$ are the zeros of $p'$? Suppose more
generally that $f$ is a polynomial with only real zeros. There is a
common misconception, reinforced by looking at small degree examples
and badly drawn graphs, that the zeros of $f'$ should be located
close to the midpoint between neighboring zeros of~$f$.  A more
accurate guide is the following:

\begin{rprinciple}  Zeros of $f'$ try to move away from the zeros of~$f$.
\end{rprinciple}

\begin{figure}[htp]
\begin{center}
\scalebox{1.3}[1.3]{\includegraphics{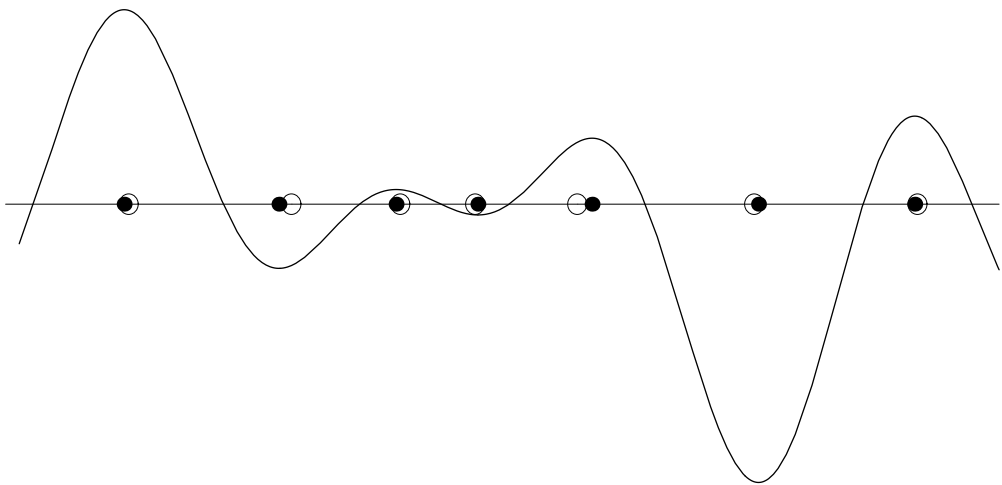}}
\caption{\sf The circles are the midpoints of neighboring zeros of the polynomial~$f$,
and the dots are the zeros of~$f'$.  The zeros of $f$
push the zeros of $f'$ away from the midpoints.} \label{fig:zfunction}
\end{center}
\end{figure}

What the principle says is:  when $f'$ decides where to put a
zero, it looks around to see what zeros of $f$ are nearby.
If it sees a bunch of zeros of $f$ to the right, it will try to shift the zero
of~$f'$ to the left, and vice-versa.
This is illustrated in Figure~\ref{fig:zfunction}, which shows
part of the graph of a high-degree polynomial.

Let's focus on the interval $(0,1)$ for the polynomial $p(x)$ given above.
Since there are many zeros to
the right, and only one zero to the left, the Repulsion Principle
says that $p'$ will shift the zero to the left, so the zero will be
closer to 0 than to 1. The question is: how close?  We answer that
question in the next section.

The Repulsion Principle, not by that name, is discussed extensively
in our paper~\cite{FR}.  This principle explains, for example,
the observation in~\cite{Mel} that the smallest zero of $f'$
is closer to the smallest zero of $f$ than to the next-smallest
zero.   The specific case of the polynomial~$p_N$
is discussed in~\cite{R}.

\section{Finding zeros of $p'$}\label{sec:logDeriv}
Any attempt to find a simple formula for the location of the first
zero of $p'$ will fail for two reasons.  First, there is no simple
formula for the derivative~$p'$.  We discuss this difficulty in the
next section.  Second, there is no closed form expression for the
zeros of an arbitrary high degree polynomial, so even if we had a
nice formula for $p'(x)$ it might not do us any good.

\subsection{A formula for $p'(x)$}

To see the difficulty of finding a formula for $p'(x)$, first consider
the product representation~\eqref{eqn:product}.  To differentiate
a product of $N+1$ things, you have to use the product rule $N$ times.
The result will be a sum of $N+1$ terms, where each term is
an $N$-fold product.  The entire expression is quite unwieldy.

Since differentiating a sum is easier than differentiating a product,
we can try multiplying out~\eqref{eqn:product} to write $p(x)$ as a sum:
\begin{equation}\label{eqn:sumform}
p(x)=\sum_{n=0}^{N+1} a_n x^n.
\end{equation}
A few of the coefficients are easy to find.  For example, $a_0=0$
because $p(0)=0$, which we see by substituting $x=0$ in the product.
It is not too hard to see that~$a_{1}=(-1)^N N!$ and $a_{N+1}=1$. A
couple more coefficients can be described as simple formulas, but
there is no known closed formula for all of the coefficients. These
coefficients come up in a variety of situations and are know as the
\emph{signed Stirling numbers of the first kind}.  We have
$a_n=s(N+1,n)$  in the usual notation for the Stirling numbers.
Since, there is no simple formula for $s(N+1,n)$,
differentiating~\eqref{eqn:sumform} does not give a useful
expression. For more information about the Stirling numbers, see
Chapter~6 of~\cite{GKP}.

Since the product form of $p$ is difficult to differentiate,
and the sum form is easy to differentiate but has unusable coefficients,
what can we do?  The answer is:
take logarithms!
Since the $\log$ of a product is a sum, we have
\begin{eqnarray}
\log(p(x))&=&\log\(\prod_{n=0}^N (x-n) \)\cr
&=& \sum_{n=0}^N \log(x-n).
\end{eqnarray}
Note: here ``$\log$'' means natural log, the function that calculator
manufacturers and nonmathematicians call ``$\ln$''.

Differentiating the above expression we
obtain ``the logarithmic derivative'' of~$p$:
\begin{equation}\label{eqn:logderivp}
\frac{p'(x)}{p(x)}= \sum_{n=0}^N \frac{1}{x-n}.
\end{equation}
This is useful because
the solutions to $p'(x)=0$ are the
same as the solutions to
\begin{equation}\label{eqn:logderivpeq0}
\frac{p'(x)}{p(x)}=0,
\end{equation}
except possibly at the points where $p(x)=0$.
However, we are trying to understand the solution
to $p'(x)=0$ which is between $0$ and~$1$, so we have the
following equivalent problem:  find the solution to
\begin{equation}
\sum_{n=0}^N \frac{1}{x-n} = 0
\end{equation}
which satisfies $0<x<1$.  We will not find an exact expression
for that solution, but we will find a good approximation.

\subsection{First approximation to the solution}

In this section we work informally and make things more precise
in the following sections

Let $\alpha$ denote the solution to $p'(x)=0$ with
$0<\alpha<1$.  As described in the previous section, we must solve
\begin{equation}\label{eqn:alpha}
0 =  \frac{1}{\alpha} +
\frac{1}{\alpha - 1} + \dots + \frac{1}{\alpha - N}.
\end{equation}
Since solving~\eqref{eqn:alpha} is the same as finding a root of an $N$th
degree polynomial, we should not expect to find an exact formula for
$\alpha$, except maybe for small~$N$.
Instead, we will look for an approximation.

Recall that the Repulsion
Principle told us that $\alpha$ is close to zero,
so in~\eqref{eqn:alpha} the term $1/\alpha$
is much larger than the rest.
So we will rearrange \eqref{eqn:alpha} to put $1/\alpha$ by itself:
\begin{equation}\label{eqn:oneoveralpha}
\frac{1}{\alpha} =
\frac{1}{1-\alpha} +\frac{1}{2-\alpha} + \dots + \frac{1}{N-\alpha}.
\end{equation}
Since $\alpha$ is close to $0$,
the right side of~\eqref{eqn:oneoveralpha}
will be close to the \emph{$N$th harmonic number}:
\begin{equation}\label{eqn:harmonicN}
H_N := 1+\frac{1}{2} + \dots + \frac{1}{N}.
\end{equation}
We give more information about the harmonic numbers in
the next section, but for now just recall that
$H_N$ is approximately $\log N$.
Thus,
\begin{equation}
\frac{1}{\alpha} \approx \log N .
\end{equation}
so, $\alpha\approx 1/\log N$, and we see that, for large $N$, $\alpha$ is actually
quite close to~$0$. In particular, $\alpha$ is not close to the
midpoint between the first two zeros.

The above calculation is informal and we do not intend to convey a precise
meaning with the ``approximately'' symbol ``$\approx$''. We
will do the calculation rigorously in the next section.

\section{Better approximations to the first zero of $p'$}

We make precise the argument in the previous section, and
find that it touches on many combinatorial and
number-theoretic objects that have been widely studied.
In our approximation to the first root of
$p'(x)$ we will encounter Euler's constant $\gamma$ and
the Riemann zeta function
\begin{equation}
\zeta(s)=\sum_{k=1}^\infty \frac{1}{k^s}.
\end{equation}
Then in Section~\ref{sec:zeta(2)} we will differentiate another
polynomial in order to evaluate $\zeta(2)$, and in Section~\ref{sec:zeta(2n)}
we evaluate the zeta function at the even integers.  That calculation
will involve the Bernoulli numbers, the Stirling numbers, and the Eulerian
numbers.

To give a precise estimate for $0<\alpha<1$, the first root of $p'$,
we need two facts.
First, the geometric series
\begin{equation}\label{eqn:geometric}
\frac{1}{1-x}
=
1+x+x^2+\cdots + x^M + O(x^{M+1}),
\end{equation}
which can be manipulated to give
\begin{equation}\label{eqn:jgeometric}
\frac{1}{j-x}
=
\frac{1}{j}+\frac{x}{j^2}+\frac{x^2}{j^3}+\cdots
    + \frac{x^M}{j^{M+1}} + O\(\frac{x^{M+1}}{j^{M+2} }\),
\end{equation}
as $x\to 0$.
Second,
\begin{equation}\label{eqn:harmonic}
H_N:=\sum_{n=1}^N \frac{1}{n}
=
\log(N) + \gamma +
         O\(\frac{1}{N}\) .
\end{equation}
Here $\gamma$ is Euler's constant, $\gamma\approx 0.577$.
Equation~\eqref{eqn:harmonic} is the \emph{definition} of Euler's
constant, but it appears in many other places, such as the Laurent
series for the Riemann-zeta function:
\begin{equation}
\zeta(s) = \frac{1}{s-1}+\gamma + O(s-1).
\end{equation}
There is a whole book on this interesting number~\cite{H}.

We have
\begin{theorem}\label{alphaApprox}
If $0<\alpha<1$ is a root of $p'_N(x)$ then
\begin{equation}
\frac{1}{\alpha} = \log(N) + O(1),
\end{equation}
so
\begin{equation}
\alpha = \frac{1}{\log(N)} + O\(\frac{1}{\log^2(N)}\).
\end{equation}
\end{theorem}
\begin{proof}
The only fact we will use about $\alpha$ is that $0<\alpha<1$,
so in particular $\alpha=O(1)$.
Start with equation \eqref{eqn:oneoveralpha},
and use equation~\eqref{eqn:jgeometric} with $M=0$, and~\eqref{eqn:harmonic}:
\begin{align*}
\frac{1}{\alpha} =&\, \sum_{j=1}^N \frac{1}{j -\alpha} \\
=&\, \sum_{j=1}^N\( \frac{1}{j} + O\(\frac{\alpha}{j^2}\)\)\\
=&\,\sum_{j=1}^N \frac{1}{j} + O(\alpha) O\(\sum_{j=1}^N
\frac{1}{j^2}\),
\\
=&\, \log(N) + \gamma+O\(\frac{1}{N}\) + O(\alpha),
        \ \ \ \ {\text{ since $\sum_n 1/n^2$ converges}}\\
=&\, \log(N) + O(1).
\end{align*}
That gives the first equation in the theorem.  To obtain the
second, take the reciprocal of both sides and use
\eqref{eqn:jgeometric} with $j=\log N$ and $x=O(1)$.
\end{proof}

We started with the assumption $0<\alpha<1$ and found an
approximation for $\alpha$ with an error term $O(1/\log^2(N))$.  By
feeding this into an almost identical calculation, we can obtain an
even better approximation to~$\alpha$.

\begin{theorem}\label{alphaApprox2}
If $0<\alpha<1$ is a root of $p'_N(x)$ then
\begin{equation}
\frac{1}{\alpha} = \log(N) + \gamma + \frac{\zeta(2)}{\log(N)} +
O\(\frac{1}{\log^2(N)}\),
\end{equation}
so
\begin{equation}
\alpha = \frac{1}{\log(N)} - \frac{\gamma}{\log^2(N)} +
\frac{\gamma- \zeta(2)}{\log^3(N)} + O\(\frac{1}{\log^4(N)} \).
\end{equation}
\end{theorem}

Euler proved that
\begin{equation}
\zeta(2) := \sum_{n=1}^\infty \frac{1}{n^2} = \frac{\pi^2}{6} .
\end{equation}
In the proof below we will only use the fact that
\begin{equation}
\sum_{n=1}^N \frac{1}{n^2} = \sum_{n=1}^\infty \frac{1}{n^2} -
\sum_{n=N+1}^\infty \frac{1}{n^2} = \zeta(2)+ O\(\frac{1}{N}\) ,
\end{equation}
which follows from \begin{equation}\frac{1}{N+1} = \int_{N+1}^\infty
\frac{1}{x^2}dx < \sum_{n=N+1}^\infty \frac{1}{n^2} <
\int_{N}^\infty \frac{1}{x^2}dx = \frac{1}{N}.\end{equation}

\begin{proof}[Proof of Theorem \ref{alphaApprox2}]
Again we start with equation \eqref{eqn:oneoveralpha}
and use equation~\eqref{eqn:jgeometric}, this time with $M=1$,
and~\eqref{eqn:harmonic}.
We have
\begin{align*}
\frac{1}{\alpha}=& \sum_{j=1}^N \frac{1}{j-\alpha} \\
=& \sum_{j=1}^N\( \frac{1}{j} +
\frac{\alpha}{j^2}+O\(\frac{\alpha^2}{j^3}\)\)\\
=& \sum_{j=1}^N  \frac{1}{j} + \alpha \sum_{j=1}^N\frac{1}{j^2}+
O(\alpha^2) \sum_{j=1}^N O\(\frac{1}{j^3}\) \\
=&\log(N)+\gamma + O\(\frac{1}{N}\) +
\alpha\(\zeta(2)+O\(\frac{1}{N}\) \) +O(\alpha^2),
\end{align*}
The last step used the fact that $\sum 1/n^3$ converges.
On the final line above, replace $\alpha$ by
$1/\log N+O(1/\log^2 N)$ and $O(\alpha^2)$ by
$O(1/\log^2 N)$, both of which follow from Theorem~\ref{alphaApprox}.
Since $1/N = O(1/\log N)$ this establishes the first formula in
Theorem~\ref{alphaApprox2}.  The second formula follows in
exactly the same way as the second formula in Theorem~\ref{alphaApprox}.
\end{proof}

It is clear that we could use the result in Theorem~\ref{alphaApprox2},
feed it back into equation~\eqref{eqn:oneoveralpha}, and obtain
an even more precise formula for~$\alpha$.  Repeating this
process
would establish that
\begin{equation}
\alpha = \frac{1}{\log(N)} +\frac{c_2}{\log^2(N)} + \cdots +
\frac{c_M}{\log^M(N)}+ O_M\(\frac{1}{\log^{M+1}(N)} \),
\end{equation}
where the $c_j$ are polynomial expressions in $\gamma$ and
$\zeta(2),\ldots,\zeta(j-1)$.  For the first few $c_j$ we have
\begin{align*}
c_2=& -\gamma \\
c_3 =& \gamma^2 - \zeta(2) \\
c_4 =& \zeta(3) - \gamma \zeta(2) + \zeta(2) - \gamma.
\end{align*}

In the next section we use a different polynomial, but which also has
equally spaced zeros, to prove Euler's result~$\zeta(2)=\pi^2/6$.


\section{The polynomials $z^N - 1$, and $\zeta(2)$}\label{sec:zeta(2)}
The previous section we considered polynomials with equally spaced
zeros on the real line.  In this section we consider polynomials with
equally spaced zeros on the unit circle.  That is,
$q(z)=q_N(z) := z^N-1$.
If we set $e(z):= e^{2 \pi i z}$ then the zeros of $q(z)$ are
at the $N$th roots of unity:  $1$, $e(1/N)$, $e(2/N)$,\ldots,$e((N-1)/N)$.
Thus,
\begin{equation}\label{eqn:unityfactor}
z^N-1 = \prod_{j\bmod N} \(z-e(j/N)\).
\end{equation}
Here $j\bmod N$ refers to a set of representatives of
the integers modulo~$N$.  This is legitimate because $e(z)$ has period~$1$.
Later we will use the specific choice $-N/2 < j \le N/2$.

We now use the polynomial $q(z)$ to evaluate $\zeta(2)$. As in the
previous section, we begin with the logarithmic derivative of~$q$.
Since we are trying to prove an identity, in this case we will make
use of both the product and the sum representation of the
polynomial. Taking the \emph{second} derivative of the logarithm
of~\eqref{eqn:unityfactor} we obtain
\begin{equation}\label{eqn:logder2}
\frac{-N(N-1)z^{N-2}-Nz^{2N-2}}{(z^N-1)^2}=
-\sum_{j\bmod N} \frac{1}{(z-e(j/N))^2} .
\end{equation}
Various formulas can be obtained by plugging in specific quantities for
$z$ in the above equation.  It is tempting to plug in $z=1$, which
is possible if one first moves the $j\equiv 0 \mod N$ term to the left
side and simplifies.  Equivalently, start over with the polynomial
$(z^N-1)/(z-1)$.   Instead, we will plug in $z=e(1/2N)$, a $2N$th root of unity,
which is motivated by noting that the left side will become
fairly simple.

Multiplying both sides of that equation by $z^2$ and plugging
in $z=e(1/2N)$ we have
\begin{lemma} We have
\begin{align}\label{deriv2qN}
\frac{N^2-2N}{4}
=&
- e(1/N) \sum_{j\bmod N} \frac{1}{(e(1/2N)-e(j/N))^2} \cr
=& - \sum_{j\bmod N} \frac{1}{\(1-e((2j-1)/2N)\)^2} .
\end{align}
\end{lemma}

The remainder of the calculation closely parallels the treatment in the
previous section: we want to expand the denominator of the right
side of \eqref{deriv2qN} as a series.  Using the
Taylor series for $e^x$ one can check that for $0<|x|<1$ we have
\begin{equation}\label{eqn:taylor}
\frac{1}{\(1-e(x)\)^2}=
-\frac{1}{4 \pi ^2 x^2}+\frac{i}{2 \pi
   x}+\frac{5}{12}-\frac{1}{6} i \pi  x+O(x^2).
\end{equation}
The restriction on $x$ is due to the fact that $e(x)=1$ for integral~$x$.
We will only need the first term in that expansion.

Since \eqref{eqn:taylor} is only valid for $0<|x|<1$ we are forced
to make the specific choice $-N/2<j\le N/2$ for the representatives
of the integers modulo~$N$. Making this choice and then applying the
Taylor expansion~\eqref{eqn:taylor} we obtain
\begin{align*}
\frac{N^2-2N}{4}
=& - \sum_{-N/2<j\le N/2} \frac{1}{(1-e((2j-1)/2N))^2} \\
=& \sum_{-N/2<j\le N/2} \( \frac{(2N)^2}{4 \pi^2 (2j-1)^2}
        +O\(\frac{N}{|2j+1|}\) \)\\
=& \frac{N^2}{ \pi^2 } \sum_{-N/2<j\le N/2} \frac{1}{(2j-1)^2}
+ O(N\log N).
\end{align*}

Dividing both sides by $N^2$ and rearranging, we have
\begin{equation}\label{idenAfterLog}
\sum_{-N/2<j\le N/2} \frac{1}{(2j-1)^2}
= \frac{\pi^2}{4} + O\(\frac{\log N}{N}\).
\end{equation}
Letting $N\to\infty$ and replacing $j\to 1-j$ for $j\le 0$, we obtain
\begin{equation}\label{eqn:z2odd}
\sum_{j=1}^\infty \frac{1}{(2j-1)^2} = \frac{\pi^2}{8} .
\end{equation}
This is equivalent to Euler's result $\zeta(2)= \pi^2/6$ because
letting
\begin{equation}
\zeta(k)=\sum_{j=1}^\infty \frac{1}{j^k}
\ \ \ \ \ \ \ \
\text{and}
\ \ \ \ \ \ \ \
\zeta_{odd}(k)=\sum_{j=1}^\infty \frac{1}{(2j-1)^k}
\end{equation}
we have the following lemma.

\begin{lemma}\label{zeta_odd} For $k>1$,
\begin{equation}
\zeta(k) = \frac{2^k }{2^{k}-1} \,\zeta_{odd}(k).
\end{equation}
\end{lemma}
\begin{proof}
\begin{align*}
\zeta(k) =&\sum_{n=1}^\infty \frac{1}{n^k} = \sum_{n \text{ odd}}
\frac{1}{n^k} + \sum_{n\text{ even}} \frac{1}{n^k}\\
=& \zeta_{odd}(k) + \sum_{n=1}^\infty \frac{1}{(2n)^k} \\
=& \zeta_{odd}(k) + \frac{1}{2^k} \zeta(k).
\end{align*}
Now solve for $\zeta(k)$.
\end{proof}

Therefore by~\eqref{eqn:z2odd},
\begin{equation}
\zeta(2) = \frac{4}{3}\zeta_{odd}(2) = \frac{\pi^2}{6}.
\end{equation}

Note that substituting $z=e(t/N)$, $0<t<1$,  into~\eqref{eqn:logder2} and 
applying the same manipulations as above gives
\begin{equation}
\sum_{n\in \Z} \frac{1}{(n+t)^2} = \frac{\pi^2}{\sin^2(\pi t)} .
\end{equation}
From this one can produce results analogous to Euler's, but where the
sum is over an arithmetic progresion.

In the next section we generalize the method to evaluate the
zeta-function at the positive even integers.

\section{Zeta at the Even Integers}\label{sec:zeta(2n)}
The same calculations as above with a little more book keeping
allow us to compute $\zeta(2n)$ for all positive integers $n$.  This general
result was known to Euler and there are many proofs of it.  To state
the result we need to introduce the Bernoulli numbers. The Bernoulli
numbers $B_k$ for $k\ge 1$ are defined as the coefficients of the
Taylor series for \begin{equation}\frac{z}{e^z-1} = \sum_{k\ge 1} \frac{B_k}{k!}
z^k =
1-\frac{1}{2}z+\frac{1}{6}\frac{z^2}{2!}-\frac{1}{30}\frac{z^3}{3!}+\cdots
.\end{equation}

For more about the Bernoulli numbers, see~\cite{GKP}.
We can now state
the general theorem which gives the values for $\zeta(2n)$.
\begin{theorem}\label{thm:zeta(2n)}
For $n\ge 1$,
\begin{equation}
\zeta(2n) := \sum_{k\ge 1}\frac{1}{k^{2n}} = \frac{(-1)^{n+1}2^{2n-1}
\pi^{2n} B_{2n}}{(2n)!}.
\end{equation}
\end{theorem}

We provide a proof of Theorem~\ref{thm:zeta(2n)} by
differentiating the polynomials $z^N-1$ in two different
ways.  In each case we will have an expression involving recursively
defined constants that are well known combinatorial objects.  Using
a couple of identities for these numbers and letting $N$ tend to
infinity we will be able to evaluate $\zeta(2n)$.

\subsection{ Repeatedly ``Differentiating''}
As in the calculation of $\zeta(2)$,
we will see that taking the logarithm of our polynomial before
differentiating will allow us to express things in a  useful
way.  However, unlike that calculation we will not simply take the
derivative. Instead of repeatedly applying the operator
$\frac{d}{dz}$ we apply the operator $z \frac{d}{dz}$.
We will refer to this as ``differentiation,'' where the quotes
indicate that we are applying $z \frac{d}{dz}$
instead of $\frac{d}{dz}$.

We wish to apply $z\frac{d}{dz}$ many times to the function
\begin{equation}\label{log(z^N-1)}
\log(z^N-1) = \sum_{j \bmod{N}}
\log\(z-e\(\frac{j}{N}\)\).\end{equation}
To
understand high ``derivatives'' of the right hand side
we will consider the high ``derivatives'' of the function
$\log(z-a),$ where $a$ is a complex number.  We have:

\begin{proposition}\label{diff RHS}
Define $b_{k,m}$ by $b_{1,1} = 1$, and $b_{0,m}=0$ for all $m\ge 0$,
and \begin{equation}
b_{k,m+1} = kb_{k,m} - (k-1)b_{k-1,m}
\end{equation}
for $ m\ge 1$ and $0\le k \le m$.
We have
\begin{equation}\(z\frac{d}{dz}\)^m \log \(z-a\) = \sum_{k=1}^m b_{k,m} z^k
(z-a)^{-k}.\end{equation}
\end{proposition}
\begin{proof}
Begin by noticing that \begin{equation}\(z\frac{d}{dz}\) (z^k(z-a)^{-k}) =
k\(z^k(z-a)^{-k} - z^{k+1}(z-a)^{-k-1}\)\end{equation} for $k\ge 1$.  We will
now establish the result by induction.  It is true for $m=1$. Assuming
it is true for for all $r<m$, with $m\ge 2$, we have
\begin{align*}
\(z\frac{d}{dz}\)^m \sum_{j \bmod N} \log(z-a) =&
\(z\frac{d}{dz}\)\sum_{k=1}^{m-1} b_{k,m-1} z^k (z-a)^{-k}\\
=& \sum_{k=1}^{m-1} b_{k,m-1}\( k z^k (z-a)^{-k} -k
z^{k+1}(z-a)^{-(k+1)}\) \\
=& \sum_{k=1}^m \(kb_{k,m-1} - (k-1)b_{k-1,m-1}\)z^k (z-a)^{-k} \\
=& \sum_{k=1}^m b_{k,m}z^k(z-a)^{-k}.
\end{align*}
The first step follows by the induction hypothesis and the
final step follows by the definition of the $b_{k,m}$.
\end{proof}

Using this result we see that
\begin{align}\label{eqn:diff1}
\(z\frac{d}{dz}\)^m \sum_{j \bmod N} \log \(z-e(j/N)\) = & \sum_{j\bmod N}
\sum_{k=1}^m b_{k,m} z^k (z-e(j/N))^{-k} \nonumber \\
=& \sum_{k=1}^m b_{k,m} \sum_{j\bmod N} z^k (z-e(j/N))^{-k},
\end{align}
which is an expression for the ``derivatives'' of the right
side of~\eqref{log(z^N-1)}.
Next we will find an expression for the
derivatives of the left side.

\begin{proposition}\label{diff LHS} Define $c_{k,m}$ by
$c_{1,1} = 1$, $c_{k,1}=0$ for all $k>1$,  and  $c_{0,m} =0$ for all
$m\ge 0$, and
\begin{equation} c_{k,m+1} =-\( k c_{k,m} + (m+1 - k) c_{k-1,m}\)
\end{equation}
for $k\ge 1$ and  $m\ge 1$. We have
\begin{equation}
\(z\frac{d}{dz}\)^m \log \(z^N-1\) = \frac{N^m \sum_{k=1}^m c_{k,m}
z^{Nk}}{(z^N-1)^m}.
\end{equation}
\end{proposition}

\begin{proof}
Define $Q_m(z) := (z^N-1)^m\(z\frac{d}{dz}\)^m \log \(z^N-1\).$  An
easy inductive argument shows that $Q_m(z)$ is a polynomial. In the
notation of the Proposition we wish to prove that \begin{equation}Q_m(z) = N^m
\sum_{k=1}^m c_{k,m} z^{Nk}.\end{equation} We will obtain a recursion for $Q_m$
and then using that recursion we will find a recursion for its
coefficients, $c_{k,m}$.

By definition we have \begin{align*}\frac{Q_{m+1}(z)}{(z^N-1)^{m+1}}
=& \(z\frac{d}{dz}\) \frac{Q_{m}(z)}{(z^N-1)^{m}} \\
=& \frac{z(z^N-1) Q_m'(z) - mNz^N Q_m(z)}{(z^N-1)^{m+1}}
.\end{align*} Thus we obtain the recursion \begin{align*}Q_{m+1}(z)
=& z(z^N-1) Q_m'(z) - mNz^N Q_m(z) \\
=& z^{N+1}Q_m'(z) - zQ_m'(z) - mNz^N Q_m(z).\end{align*}

By induction we have the following
\begin{align*}
z^{N+1} Q_m'(z) =&  N^{m+1} \sum_{k=1}^m k c_{k,m} z^{N(k+1)} \\
zQ_m'(z)=& N^{m+1} \sum_{k=1}^m kc_{k,m}z^{Nk}\\
mNz^N Q_m(z) =& N^{m+1} \sum_{k=1}^m m c_{k,m} z^{N(k+1)}
\end{align*}
Combining these we obtain the recursive formula for the constants,
$c_{k,m}$.
\end{proof}

\subsection{Combinatorial Numbers}
In Propositions~\ref{diff RHS} and~\ref{diff LHS} we introduced two
sets of numbers $\{b_{k,m}\}$ and $\{c_{k,m}\}$, and those numbers
satisfied particular recursion formulas. It turns out that those
numbers are related to other famous numbers from combinatorics:
Stirling numbers of the second kind and the Eulerian numbers.

The Stirling numbers of the second kind $S(n,k)$ count the number of
ways to partition a set of $n$ elements into $k$ nonempty subsets.
They can be described recursively by $S(1,1) = 1$ and the recursion
\begin{equation}
S(m,k) = S(m-1,k-1)+kS(m-1,k).
\end{equation}
On the other hand, the Eulerian number $e_{k,m}$ is the number of
permutations of the numbers $1$ to $m$ in which exactly $k+1$
elements are greater than the previous element.  The Eulerian
numbers are given recursively by $e_{1,1}=1$ and $e_{k,1} = 0$ for
$k> 1$, and
\begin{equation}
e_{k,m} = (m-k+1)e_{k-1,m-1} + k e_{k,m-1}.
\end{equation}
It is the recursive definitions which we will make use of.

We have the following:
\begin{proposition}\label{combinatorics} In the notation above,
\begin{align*}b_{k,m} =& (m-1)!S(m,k)\\
c_{k,m}=& (-1)^{m-1} e_{k,m}.\end{align*} In particular, $b_{m,m} =
(-1)^{m-1}(m-1)!$.
\end{proposition}

\begin{proof}
We will not prove the first of these two results because the only
fact we need about $b_{n,m}$ is $b_{m,m} = (-1)^{m-1}(m-1)!$. To
prove this, it is enough to note that by
the recursion in Proposition~\ref{diff RHS} we have
\begin{equation}b_{m,m} = mb_{m,m-1} - (m-1)b_{m-1,m-1} = -(m-1)b_{m-1,m-1}.\end{equation}
Thus the result follows by induction and the fact that $b_{1,1} =
1$.

By the definition of $c_{k,m}$ in Proposition~\ref{diff LHS} we see
that $c_{k,m}$ and $(-1)^{m-1}e_{k,m-1}$ satisfy the same recurrence
relation with the same initial conditions.
\end{proof}

The following
relates the Eulerian numbers that appear in our formulas to the
Bernoulli numbers which appear
in the formula for $\zeta(2n)$.
This Lemma is Exercise~72 in Chapter~6 of~\cite{GKP}.
\begin{lemma}\label{eulerian_bernoulli}
\begin{equation}\sum_{k=1}^n (-1)^k c_{k,n} = (-1)^n  2^{n} (2^{n} -
1)\frac{B_{n}}{n}.\end{equation}
\end{lemma}

\subsection{Limit as $N\to \infty$}

We now return to our discussion of the two expressions for the
``derivatives'' of $\log(z^N-1)$.  As in the calculation of
$\zeta(2)$, to evaluate $\zeta(2n)$ we begin by equating the two
expressions for $\(z\frac{d}{dz}\)^{(2n)} \log(z^N-1)$ given by
equation \eqref{eqn:diff1} and Proposition~\ref{diff LHS}.  Using
the Taylor series expansion that we used in the previous section and
similar analysis we are led to the evaluation of $\zeta_{odd}(2n)$.

\begin{proof}[Proof of Theorem \ref{thm:zeta(2n)}] Set $m=2n$.
Beginning with equation \eqref{eqn:diff1} and setting $z=e(1/2N)$ we
have
\begin{align*}
\(z\frac{d}{dz}\)^m \sum_{j \bmod N} \log \(z-e(j/N)\)\mid_{z=e(1/2N)}
=&
\sum_{k=1}^m b_{k,m} \sum_{j\bmod N} e(1/2N)^k (e(1/2N)-e(j/N))^{-k} \\
=& \sum_{k=1}^m b_{k,m} \sum_{j\bmod N}
\(1-e\(\frac{(2j-1)}{2N}\)\)^{-k}
\end{align*}

As in the calculation from Section~\ref{sec:zeta(2)} we will use the
Taylor series
\begin{equation} (1-e(x))^{-1} = \frac{-1}{2\pi i x} + O(1),\end{equation} which is valid for
$0<|x|<1$. To use the Taylor series we must have $-N/2 < j \le N/2$ and
so we obtain
\begin{align}\label{eqn:1}
\sum_{k=1}^m b_{k,m} \sum_{-N/2 < j\le N/2}&
\(1-e\(\frac{(2j-1)}{2N}\)\)^{-k} \nonumber \\
=& \sum_{k=1}^m
b_{k,m} \sum_{-N/2 < j\le N/2} \left[\(\frac{2N}{2\pi i (2j-1)}\)^k
+
O\( \frac{N^{k-1}}{(2j-1)^{k-1}}\) \right] \nonumber\\
=& b_{m,m}  \sum_{-N/2 < j\le N/2} \(\frac{2N}{2\pi i (2j-1)}\)^m +
O_m\(N^{m-1}\) \nonumber \\ &+ \sum_{k=2}^{m-1} b_{k,m} \sum_{-N/2 <
j\le N/2} \left[ \(\frac{2N}{2\pi i (2j-1)}\)^k + O_k\(
\frac{N^{k-1}}{(2j-1)^{k-1}}\)\right] \nonumber \\ &+
\sum_{-N/2 < j\le N/2} \left[\(\frac{2N}{2\pi i (2j-1)}\) + O\(1\) \right]\nonumber  \\
=&b_{m,m}\(\frac{N}{\pi i}\)^m \sum_{-N/2 < j\le N/2}
\(\frac{1}{(2j-1)}\)^m + O_m\(N^{m-1}\) + O(N \log(N)).
\end{align}
We write $O_m$ because the implied constant depends on $m$.
The $O(N\log N)$ term is only necessary when~$m=2$.

Next we will turn to our other expression for
$\(z\frac{d}{dz}\)\log(z^N-1),$ which was given in Proposition~\ref{diff LHS}.  As before, evaluating at $e(1/2N)$ we have
\begin{align}\label{eqn:2}
\frac{N^m \sum_{k=1}^m c_{k,m}
e(1/2N)^{Nk}}{(e(1/2N)^N-1)^m} =& \frac{N^m}{(-2)^m} \sum_{k=1}^m
c_{k,m} (-1)^k \nonumber\\
=&\frac{N^m}{2^m} \( (-1)^m 2^{m} (2^{m} - 1)\frac{B_{m}}{m }\)
\nonumber
\\ =& N^m\((2^{m} - 1)\frac{B_{m}}{m }\),
\end{align}
where the second equality follows from the Lemma~\ref{eulerian_bernoulli} and we use that $m$ is even.
By equation~\eqref{log(z^N-1)} we have that equations~\eqref{eqn:1}
and~\eqref{eqn:2} are equal.  Setting the last lines of each
equation equal and dividing by $N^m$ we have
\begin{equation}\( (2^{m} - 1)\frac{B_{m}}{m }\) =
b_{m,m}\({\pi i}\)^{-m} \sum_{-N/2<j\le N/2} \(\frac{1}{(2j-1)}\)^m +
O_m\(N^{-1}\) + O(N^{1-m} \log(N)).
\end{equation}
Letting $N$ go to infinity we obtain
\begin{equation}\label{doubleSum}
\frac{b_{m,m}}{(\pi i)^m}\sum_{j=-\infty}^\infty \frac{1}{(2j-1)^m}
= \frac{ (2^m - 1)  B_m }{m}.
\end{equation}
To finish the proof we use the facts that
$b_{m,m} = (-1)^{m-1}(m-1)!$
from Lemma~\ref{diff RHS} and Lemma~\ref{zeta_odd} to obtain
\begin{align*} \zeta(m) =& \frac{2^m \zeta_{odd}(m) }{(1-2^m)} =
\frac{-2^{m-1} \pi^m i^m B_m}{m!} \\=& \frac{2^m
}{(1-2^m)}\frac{1}{2}\sum_{j=-\infty}^\infty \frac{1}{(2j-1)^m} \\
=& \frac{2^{m-1} }{(1-2^m)} \frac{ (2^m - 1)(\pi i)^m  B_m }{m! }.
\end{align*}
Using $m=2n$ we have
\begin{equation}\zeta(2n) = \frac{\pi^{2n} (-1)^{n-1}2^{2n-1} B_{2n}}{(2n)!},\end{equation}
the desired result.
\end{proof}

Instead of the substitution $z=e(1/2N)$ one can set $z=e(A/B N)$ to obtain an
expression for
\begin{equation}
\sum_{n\in \Z} \frac{1}{(A+B n)^k} .
\end{equation}
This can then be used to evaluate Dirichlet L-functions $L(s,\chi)$,
with character~$\chi$,
at the positive even~(odd) integers if $\chi$ is even~(odd).

\end{document}